\documentclass[11pt]{amsart}
\usepackage[T1]{fontenc}
\usepackage{lmodern}
\usepackage{amsmath,amssymb,amsthm,booktabs,longtable}
\usepackage[margin=1.1in]{geometry}
\usepackage[hidelinks]{hyperref}
\newtheorem{theorem}{Theorem}[section]
\newtheorem{lemma}[theorem]{Lemma}
\newtheorem{proposition}[theorem]{Proposition}

\newtheorem{conjecture}[theorem]{Conjecture}
\theoremstyle{definition}
\newtheorem{definition}[theorem]{Definition}
\newtheorem{remark}[theorem]{Remark}
\newcommand{\N}{\mathbb N}
\newcommand{\ord}{\operatorname{ord}}
\newcommand{\core}{\operatorname{core}}
\newcommand{\up}{\mathord{\uparrow}}
\newcommand{\eps}{\varepsilon}
\newcommand{\ind}[2]{[#1:#2]}
\title{The Herzog--Sch\"onheim conjecture for solvable groups}
\author{Murali Menon}
\date{October 1, 2026}

\begin{document}
\begin{abstract}
The Herzog--Sch\"onheim conjecture states that if a group is partitioned
into finitely many, and at least two, cosets of subgroups, then two of
these subgroups have the same index. We prove the conjecture for all
solvable groups, finite or infinite. Sun's bound on unions of cosets,
applied in a quotient, is combined with divisibility restrictions on
local indices relative to the subgroup generated by the $p$-elements,
where $p$ is the largest prime divisor of the order of a minimal
counterexample. These restrictions are obtained by counting fixed points
in a primitive solvable quotient. This suffices except at the prime $7$,
where the counting bound narrowly fails and an additional
overlap argument is needed. The proof has been formalized in Lean~4 with
mathlib.
\end{abstract}
\maketitle

\section{Introduction}\label{sec:intro}

Erd\H{o}s conjectured that in any partition of the integers into finitely
many, and at least two, arithmetic progressions, the largest common
difference occurs at least twice. This was proved independently by
Mirsky and Newman and by Davenport and Rado; see the introduction of
\cite{GM} for the history. Viewing arithmetic progressions as cosets of subgroups of
$\mathbb Z$, Herzog and Sch\"onheim proposed the following generalization
\cite{HS}.

\begin{conjecture}[Herzog--Sch\"onheim]\label{conj:HS}
Let $G$ be a group and let
\begin{equation}\label{eq:partition}
 G=\bigsqcup_{i=1}^k g_iH_i,
 \qquad H_i\leq G,\quad g_i\in G,
\end{equation}
be a partition of $G$ into $k\geq2$ left cosets. Then
$\ind G{H_i}=\ind G{H_j}$ for some $i\ne j$.
\end{conjecture}

We call the cosets $g_iH_i$ the \emph{parts} of the partition.
Following \cite{GM}, we say that $G$ is \emph{HS} if every such partition
of $G$ contains two cosets whose underlying subgroups have equal index.
The conjecture is known for finite
nilpotent groups \cite{BFF86} and pyramidal groups \cite{BFF87}, for
partitions by subnormal subgroups \cite{Sun}, under conditions on the
prime divisors of $|G|$ \cite{GS}, for groups of order less than $1440$
\cite{MS}, and for simple and symmetric groups \cite{GM}. For arbitrary
groups, Akman and Sissokho proved it for partitions into cosets of at
most seven distinct subgroups \cite{AS25a}. They also study
\emph{Steiner} partitions, which use exactly one coset of each subgroup,
when the subgroups commute pairwise; the case of five mutually commuting
subgroups is treated in \cite{AS26}; see also \cite{AS25b}. A preprint
claiming the solvable case was withdrawn \cite{Burkhart}; we do not use
its argument. Our main result is the following.

\begin{theorem}\label{thm:main}
Every solvable group is HS.
\end{theorem}

\subsection*{The idea}
For a finite group, a partition \eqref{eq:partition} satisfies
$\sum_i1/\ind G{H_i}=1$. Garonzi and Margolis prove the conjecture for
simple and symmetric groups by showing that the sum of $1/m$ over the
distinct subgroup indices $m$ of $G$ is less than $2$
\cite[Lemma~2.1]{GM}. Such a criterion depends only on the set of subgroup
indices, and it cannot work for solvable groups: the sum is unbounded
already on cyclic groups, as they note. Our argument also uses where the
cosets lie. We compare the reciprocal sum of a family of parts with the
size of the union of their images in a quotient $G/F$, and bound that
size from below by a theorem of Sun on unions of cosets.

For a finite group $G$, put
$\Pi(G)=\prod_{q\mid |G|}(1-q^{-1})^{-1}$, with $q$ ranging over primes.
Ginosar and Schnabel showed that $\Pi(G)\leq2$ implies that $G$ is HS
\cite[Theorem~A]{GS}. Our criterion (Proposition~\ref{prop:criterion})
takes $p$ to be the largest prime divisor of $|G|$. It combines the
corresponding product over the primes below $p$ with a factor that
depends on the multiplicative orders of the primes $q<p$ modulo $p$
and tends to zero with $p$. That factor comes from the local indices,
which are not visible to a criterion on $|G|$ alone.

\subsection*{How the proof works}
Let $G$ be a finite solvable group of least order that is not HS, and
let $d_i=\ind G{H_i}$. No two $d_i$ are coprime, because cosets of
subgroups of coprime index always meet. Using the minimality of $|G|$,
an elementary reduction (Lemma~\ref{lem:mincex}) also shows that no
$d_i$ equals $2$ and that the largest prime $p$ dividing $|G|$ is at
least $5$. The rest of the proof plays two facts
against each other. The indices must be spread out, since they are
distinct numbers whose reciprocals sum to $1$. But the group structure
forces them into a thin set.

\emph{Step 1: split each index at a normal subgroup}
(Sections~\ref{sec:local} and~\ref{sec:structure}). Let $F$ be the
subgroup generated by the elements of $p$-power order in $G$. It is a
nontrivial normal subgroup, and it has no proper normal subgroup of index
prime to $p$. Each index factors as
\[
 d_i=h_i\ell_i,
 \qquad h_i=\ind G{H_iF},\qquad \ell_i=\ind F{H_i\cap F}.
\]
We call $h_i$ the \emph{quotient index} and $\ell_i$ the \emph{local
index} of the part. The part meets each $F$-coset it touches in
$|F|/\ell_i$ elements. A part with $\ell_i=1$ is a union of $F$-cosets;
we call it \emph{universal}. If every part is universal, the partition
descends to the smaller group $G/F$ with the same indices, contrary to
minimality. So the set $J$ of indices of the non-universal parts is nonempty.

\emph{Step 2: Sun's bound applied in the quotient}
(Theorem~\ref{thm:sun} and Lemma~\ref{lem:shadow}). Record each positive
integer $n$ by its exponent vector $v(n)$, the list of exponents in its
prime factorization, and give that vector the weight $1/n$. Let $U$ be
the set of vectors that are coordinatewise at least $v(h_i)$ for some
$i\in J$; these are the vectors of the multiples of the quotient indices.
Write $W(U)$ for the total weight of $U$. The non-universal parts
together cover a union of $F$-cosets, so their images cover $r\nu$ of
the $r=\ind GF$ cosets in $G/F$, where $\nu=\sum_{i\in J}1/d_i$. Sun's
lower bound for the size of a union of cosets then gives
$W(U)\leq\Pi\nu$, where $\Pi$ depends only on the primes involved.

\emph{Step 3: local indices are restricted} (Lemma~\ref{lem:cyclotomic}).
For distinct primes $p,q$, write $\ord_p(q)$ for the multiplicative order
of $q$ modulo $p$. If $p\nmid\ell_i$, a Sylow $p$-subgroup of $F$ lies in a
maximal subgroup $M$ of $F$ containing $H_i\cap F$. The primitive solvable
quotient of $F$ acting on the cosets of $M$ has prime-power degree $q^c$,
and the Sylow subgroup acts on its minimal normal subgroup
nontrivially, because $F$ has no proper normal subgroup of index prime to
$p$. Counting fixed points shows that each $\ell_i>1$ is divisible by $p$
or by $q^{\ord_p(q)}$ for some prime $q<p$.

\emph{Step 4: counting gives a contradiction for $p\ne7$}
(Proposition~\ref{prop:criterion} and Lemma~\ref{lem:numeric}).
Since $d_i=h_i\ell_i$, we have
$v(d_i)=v(h_i)+v(\ell_i)$. By Step~3, every $v(d_i)$ with $i\in J$ lies
either in $U+v(p)$ or in $U+v(q^{\ord_p(q)})$ for some prime $q<p$.
Shifting a set by $v(m)$ multiplies its weight by $1/m$. Accounting for
their overlap, the union of these shifted sets has weight at most
$c_pW(U)$, where
$c_p=1/p+(1-1/p)\sum_{q<p}q^{-\ord_p(q)}$ tends to zero as $p$ grows.
The $v(d_i)$
are distinct and their weights add up to $\nu$, so
$\nu\leq c_pW(U)\leq c_p\Pi\nu$, which requires $c_p\Pi\geq1$.
Here $\Pi$ is taken over all primes at most $p$. Exact rational estimates
for the small primes, an interval argument up to $255$, and a uniform
elementary bound for $p\geq256$ show that $c_p\Pi<1$ for every prime
$p\geq5$ except $7$; the finite inequalities are also checked in the
Lean development.

\emph{Step 5: the prime $7$} (Section~\ref{sec:seven}). At $p=7$ a local
index may be divisible by $8=2^{\ord_7(2)}$, and the shifted copy
$U+v(8)$ has the comparatively large weight $W(U)/8$. As a result
$c_7\Pi$ is about $1.1$ with all four primes included, and the count of
Step~4 falls short of a
contradiction by about $10\%$. We retain information that the generic
count discards: which minimal quotient-index vectors actually carry parts
of local index exactly $7$ or exactly $8$, and which parts meet the same
$F$-coset as such a part (Lemma~\ref{lem:chambers}). A missing-point and
overlap-charging argument turns this incidence information into
incompatible bounds (Proposition~\ref{prop:seven}).

We use Sun's coset union bound and the standard structure theory of
primitive solvable groups.
The new ingredients are the combination of the local divisibility
restriction with the quotient count, and the additional argument at the
prime $7$. Section~\ref{sec:formal} describes the Lean formalization
\cite{Repository}.

\section{Preliminaries}\label{sec:prelim}

Throughout, $\N=\{0,1,2,\ldots\}$. For a subgroup $A$ of a group $G$,
its normal core is $\core_G(A)=\bigcap_{g\in G}gAg^{-1}$.

\subsection{Coset intersections and reductions}
The following standard fact is used repeatedly.

\begin{lemma}[{\cite[Lemma~2.1 and Remark~2.2]{Sun06}}]\label{lem:coprime}
Let $X$ be a finite group and $Y,Z\leq X$ with $\ind XY$ and $\ind XZ$
coprime. Then $YZ=X$, and every coset of $Y$ meets every coset of $Z$.
\end{lemma}
Sun's Lemma~2.1 states $YZ=X$ for subgroups of finite coprime index in any
group. His Remark~2.2 states that if $aY\cap bZ=\varnothing$ for some
$a,b$ then the indices are not coprime. The same statement appears in
\cite[Corollary~2.1]{GS} and, in terms of harmonic tuples, in
\cite[Lemma~2.2]{MS}.

\begin{lemma}\label{lem:finite}
If every finite solvable group is HS, then every solvable group is HS.
More precisely, for a partition \eqref{eq:partition} of a group $G$ the
subgroup $N=\core_G\bigl(\bigcap_iH_i\bigr)$ is normal of finite index,
and the cosets $g_iH_i/N$ partition $G/N$ with
$\ind{G/N}{H_i/N}=\ind G{H_i}$.
\end{lemma}
\begin{proof}
By Korec and Zn\'am \cite[Theorem~1]{KZ}, every $H_i$ has finite index
in $G$. Hence $K=\bigcap_iH_i$ has finite index, with
$\ind GK\leq\prod_i\ind G{H_i}$ \cite[Lemma~1]{KZ}, and so does its normal
core $N$, the intersection of the finitely many conjugates of $K$
\cite[Lemma~2]{KZ}. Each part $g_iH_i$ is a union of $N$-cosets, so the
images are pairwise disjoint and partition $G/N$, with the same indices.
This is \cite[Theorem~4]{KZ}. The group $G/N$ is solvable if $G$ is.
\end{proof}

\begin{lemma}\label{lem:mincex}
Let $G$ be a finite solvable group of least order that is not HS, and
fix a partition \eqref{eq:partition} with $k\geq2$ and distinct indices
$d_i=\ind G{H_i}$. Then
\begin{enumerate}
\renewcommand{\labelenumi}{(\alph{enumi})}
\item $\sum_i1/d_i=1$ and every $d_i>1$;
\item no two of the $d_i$ are coprime;
\item no $d_i$ equals $2$;
\item the largest prime divisor $p$ of $|G|$ is at least $5$.
\end{enumerate}
\end{lemma}
\begin{proof}
(a) The parts have sizes $|G|/d_i$ and are disjoint. A part of index $1$
would be all of $G$, but $k\geq2$. (The sum is \cite[Theorem~2]{KZ}.)

(b) Distinct parts are disjoint, so this follows from
Lemma~\ref{lem:coprime}. (This is also \cite[Theorem~3]{KZ}.)

(c) Suppose $d_1=2$. The other parts partition the complementary coset
$xH_1$. Each $g_jH_j\subseteq xH_1$ forces $H_j\leq H_1$, and translating
by $x^{-1}$ gives a partition of $H_1$ into cosets of the $H_j$,
$j\geq2$, with distinct indices $\ind{H_1}{H_j}=d_j/2$. If $k\geq3$ this
partition has at least two parts, contradicting the minimality of $|G|$.
If $k=2$ then $g_2H_2=xH_1$ and $d_2=2=d_1$, contradicting distinctness.
This is the reduction of \cite[Lemma~2.3]{GS}; see also
\cite[Lemma~2.3]{MS}.

(d) If $p\leq3$, every index is of the form $2^a3^b$. By (b), either all
indices are even or all are divisible by $3$: if an odd index occurs, it
is a nontrivial power of $3$ and every other index is divisible by $3$.
By (c), the index $2$ is excluded. We have the following reciprocal sums:
\[
 \sum_{\substack{a\geq1,\ b\geq0\\(a,b)\ne(1,0)}}2^{-a}3^{-b}
 =\frac1{1-1/3}-\frac12=1,\qquad
 \sum_{\substack{a\geq0,\ b\geq1}}2^{-a}3^{-b}
 =\frac1{1-1/2}\cdot\frac{1/3}{1-1/3}=1.
\]
The $d_i$ form a finite subset of one of these two infinite index sets,
so $\sum_i1/d_i<1$, contrary to (a).

Part (d) also follows from \cite[Theorem~B]{GS}: every finite group whose
order has exactly two prime divisors is HS (a $p$-group is nilpotent, so
covered by \cite{BFF86}). The argument above keeps the
proof self-contained.
\end{proof}

\subsection{Sun's coset union bound}
A set $A\subseteq\mathbb Z$ that is invariant under translation by some
positive integer $m$ has density
$\operatorname{dens}(A)=m^{-1}\#\bigl(A\cap\{0,\ldots,m-1\}\bigr)$, which
does not depend on the choice of $m$. For example,
$\operatorname{dens}(h\mathbb Z)=1/h$.

\begin{theorem}[Sun's coset union bound]\label{thm:sun}
Let $L$ be a finite solvable group and let $x_jL_j$ be a nonempty finite
family of cosets of subgroups of $L$. Then
\[
 |L|\cdot\operatorname{dens}\Bigl(\bigcup_j\ind L{L_j}\mathbb Z\Bigr)
 \leq\left|\bigcup_j x_jL_j\right|.
\]
\end{theorem}

This is the case $H=1$ of \cite[Theorem~3.1]{Sun}. Sun's theorem
concerns a subgroup $H$ of finite index and subgroups $G_1,\ldots,G_k$
containing $H$, under either of two hypotheses; the second is a
composition series from $H$ to the whole group with factors of prime
order, which a finite solvable group has. With $H=1$, Sun's conclusion
says that the union of the cosets has at least as many elements as there
are $n\in\{0,\ldots,|L|-1\}$ divisible by some $\ind L{L_j}$, and this
number is the left side above. Sun's theorem has no hypothesis that the
cosets cover the group, and we apply it to families whose union is a
proper subset. In the nilpotent case the bound goes back to Berger,
Felzenbaum and Fraenkel \cite{BFF86}.

\subsection{Exponent vectors and weights}
Fix a finite set $\mathcal P$ of
primes. For a positive integer $n$ whose prime divisors lie in
$\mathcal P$, let $v(n)=(v_q(n))_{q\in\mathcal P}$ be its exponent
vector. Order $\N^{\mathcal P}$ coordinatewise. For $a\in\N^{\mathcal P}$
and $E\subseteq\N^{\mathcal P}$, put
\[
 \up a=\{b\in\N^{\mathcal P}:a\leq b\},\qquad
 \up E=\bigcup_{a\in E}\up a.
\]
Thus $\up v(n)$ is the set of vectors of the multiples of $n$ whose
prime divisors lie in $\mathcal P$.
For $q\in\mathcal P$ let $e_q=v(q)$ be the unit vector
at $q$. A set $E\subseteq\N^{\mathcal P}$ is \emph{upward closed} if
$a\in E$ and $a\leq b$ imply $b\in E$; every set $\up E$ is upward
closed. Define
\[
 w(a)=\prod_{q\in\mathcal P}q^{-a_q},\qquad
 W(E)=\sum_{a\in E}w(a),\qquad
 \Pi=\prod_{q\in\mathcal P}(1-q^{-1})^{-1},
\]
so that $w(v(n))=1/n$ and $W(\N^{\mathcal P})=\Pi$. For a set $E$ and a
vector $b$ write $E+b=\{a+b:a\in E\}$. Then
\begin{equation}\label{eq:translate}
 W(E+b)=w(b)W(E),\qquad W(\up a)=\Pi w(a).
\end{equation}

\begin{lemma}\label{lem:dictionary}
Let $h_1,\ldots,h_s$ be positive integers whose prime divisors lie in
$\mathcal P$, and put $U=\bigcup_j\up v(h_j)$. Then
\begin{equation}\label{eq:dictionary}
 W(U)=\Pi\cdot\operatorname{dens}\Bigl(\bigcup_jh_j\mathbb Z\Bigr).
\end{equation}
\end{lemma}
\begin{proof}
This is a reformulation of \cite[Lemma~3.4]{Sun}, which Sun attributes to
an earlier paper of his; we give the short proof for completeness.
For a nonempty $S\subseteq\{1,\ldots,s\}$ let $h_S=\operatorname{lcm}_{j\in S}h_j$.
Then $\bigcap_{j\in S}\up v(h_j)=\up v(h_S)$ and
$\bigcap_{j\in S}h_j\mathbb Z=h_S\mathbb Z$. By inclusion--exclusion
and \eqref{eq:translate}, the two sides of \eqref{eq:dictionary} are
$\sum_{S\ne\varnothing}(-1)^{|S|+1}\Pi/h_S$ and
$\Pi\sum_{S\ne\varnothing}(-1)^{|S|+1}/h_S$.
\end{proof}

\section{Quotient and local indices}\label{sec:local}

In this section $G$ is a finite solvable group with a partition
\eqref{eq:partition}, $F$ is a normal subgroup of $G$, and $r=\ind GF$.
We take $\mathcal P$ to contain every prime divisor of $|G|$.

\begin{definition}\label{def:local}
For each part $g_iH_i$ put
\[
 d_i=\ind G{H_i},\qquad h_i=\ind G{H_iF},\qquad \ell_i=\ind F{H_i\cap F},
\]
so that $d_i=h_i\ell_i$. We call $h_i$ the \emph{quotient index} and
$\ell_i$ the \emph{local index} of the part. The part is
\emph{universal} if $\ell_i=1$, that is, if $F\leq H_i$. Let
$J=\{i:\ell_i>1\}$ be the set of indices of the non-universal parts, let
$X_i\subseteq G/F$ be the image of $g_iH_i$, and put
\[
 \nu=\sum_{i\in J}\frac1{d_i},\qquad
 U=\bigcup_{i\in J}\up v(h_i).
\]
\end{definition}

The image $X_i$ is a coset of $H_iF/F$ and has $r/h_i$ elements. If
$g_iH_i$ meets an $F$-coset $yF$ in a point $y'$, then
$g_iH_i\cap yF=y'(H_i\cap F)$ has $|F|/\ell_i$ elements.

\begin{lemma}\label{lem:shadow}
For every $I\subseteq J$,
\[
 W\Bigl(\bigcup_{i\in I}\up v(h_i)\Bigr)\leq
 \frac{\Pi}{r}\,\Bigl|\bigcup_{i\in I}X_i\Bigr|.
\]
Moreover $\bigl|\bigcup_{i\in J}X_i\bigr|=r\nu$, so
\begin{equation}\label{eq:shadow}
 W(U)\leq\Pi\nu.
\end{equation}
\end{lemma}
\begin{proof}
The first inequality is trivial for $I=\varnothing$. Otherwise apply
Theorem~\ref{thm:sun} in the solvable group $G/F$ to the cosets $X_i$,
$i\in I$, whose subgroups $H_iF/F$ have index $h_i$, and then
Lemma~\ref{lem:dictionary}.

An $F$-coset met by a non-universal part cannot meet a universal part,
since a universal part contains every $F$-coset it meets. Hence the
union of the non-universal parts is exactly the inverse image of
$\bigcup_{i\in J}X_i$. Its size is $|G|\nu$, which gives
$\bigl|\bigcup_{i\in J}X_i\bigr|=r\nu$.
\end{proof}

\section{The subgroup \texorpdfstring{$F$}{F} and local divisibility}\label{sec:structure}

An element of a finite group is a \emph{$p$-element} if its order is a
power of $p$. The subgroup generated by the $p$-elements is also
generated by the Sylow $p$-subgroups. For primes $p\ne q$,
$\ord_p(q)$ denotes the multiplicative order of $q$ modulo $p$.

\begin{lemma}\label{lem:structure}
Let $G$ be a finite group and $p$ a prime dividing $|G|$, and let $F$ be
the subgroup generated by the $p$-elements of $G$. Then $F$ is a
nontrivial normal subgroup of $G$, and every normal subgroup of $F$ of
index prime to $p$ is $F$ itself.
\end{lemma}
\begin{proof}
$F$ is characteristic in $G$, and it contains a Sylow $p$-subgroup of $G$,
which is nontrivial. Let $N\lhd F$ with $p\nmid\ind FN$. A $p$-subgroup $T$
of $F$ maps to a $p$-subgroup of $F/N$, whose order is prime to $p$.
Thus $T\leq N$. Each Sylow $p$-subgroup of $G$ consists of
$p$-elements, so it is a $p$-subgroup of $F$ and lies in $N$. These
subgroups generate $F$; hence $N=F$.
\end{proof}

Lemma~\ref{lem:structure} does not require solvability. Only its property
concerning normal subgroups of index prime to $p$
is used in the next lemma.

\begin{lemma}[Cyclotomic divisibility]\label{lem:cyclotomic}
Let $F$ be a finite solvable group and $p$ a prime such that every normal
subgroup of $F$ of index prime to $p$ is $F$ itself. If $Y<F$ and
$p\nmid\ind FY$, then some prime $q\ne p$ dividing $|F|$ satisfies
\[
 q^{\ord_p(q)}\mid\ind FY.
\]
\end{lemma}
\begin{proof}
Let $P$ be a Sylow $p$-subgroup of $F$. The $P$-orbits on $F/Y$ have
$p$-power sizes summing to $\ind FY$, so $P$ fixes some coset $xY$, that
is, $P\leq xYx^{-1}$. Replacing $Y$ by this conjugate does not change the
index, so we may assume $P\leq Y$. Choose a maximal subgroup $M$ of $F$
containing $Y$, let $R=\core_F(M)$, and let $V/R$ be a minimal normal
subgroup of the solvable group $F/R$. Then $V/R$ is elementary abelian of
order $q^c$ for a prime $q$, and $c\geq1$.

We have $V\nleq M$, since otherwise the normal subgroup $V$ would lie in
$R$. Hence $VM=F$ by maximality. The subgroup $V\cap M$ contains $R$, so it
is normal in $V$, as $V/R$ is abelian; it is also normalized by $M$, hence
normal in $VM=F$. As $R$ is the largest normal subgroup of $F$ contained
in $M$, this gives $V\cap M\leq R$, and therefore $V\cap M=R$. So
$\ind FM=\ind V{V\cap M}=q^c$. As $\ind FM$ divides $\ind FY$, which is
prime to $p$, we get $q\ne p$; and $q$ divides $|F|$.

Since $VM=F$ and $V\cap M=R$, the map
$\varphi:V/R\to F/M$, $vR\mapsto vM$, is a bijection.
For $x\in P\leq M$ it satisfies
\[
 \varphi(xvx^{-1}R)=(xvx^{-1})M=xvM=x\varphi(vR).
\]
Thus it carries conjugation by $P$ on $V/R$ to left multiplication by
$P$ on $F/M$. If $P$ fixed every point of $F/M$, then $P\leq R$;
but $R\lhd F$ contains a Sylow $p$-subgroup, so
$\ind FR$ is prime to $p$ and $R=F$ by hypothesis, contrary to $R\leq M<F$.
Hence the fixed subgroup $C_{V/R}(P)$ has order $q^b$ with $b<c$. The
non-singleton $P$-orbits on $V/R$ have size divisible by $p$, so
$q^c\equiv q^b\pmod p$, hence $q^{c-b}\equiv1\pmod p$ and
$\ord_p(q)\leq c-b\leq c$. Therefore
$q^{\ord_p(q)}\mid q^c=\ind FM$, which divides $\ind FY$.
\end{proof}

The argument is the usual one for solvable primitive groups, where a
minimal normal subgroup of the primitive quotient $F/R$ acts regularly;
compare \cite[Theorem~7.2.6]{Robinson}. Counting fixed points avoids a
separate appeal to the orders of general linear groups.

\section{Largest prime other than \texorpdfstring{$7$}{7}}\label{sec:general}

For a partition and normal subgroup $F$, we use the notation of
Definition~\ref{def:local}: $d_i=h_i\ell_i$, $J=\{i:\ell_i>1\}$,
$\nu=\sum_{i\in J}1/d_i$, and $U=\bigcup_{i\in J}\up v(h_i)$.
The weight $W$ and the product $\Pi$ are as in Section~\ref{sec:prelim},
for the specified prime set $\mathcal P$.

For a prime $p$ put
\[
 \Pi_p=\prod_{q<p}(1-q^{-1})^{-1},\qquad
 S_p=\sum_{q<p}q^{-\ord_p(q)},\qquad
 c_p=\frac1p+\left(1-\frac1p\right)S_p,
\]
where here and below sums and products indexed by $q$ run over primes.

\begin{proposition}\label{prop:criterion}
Let $G$ be a finite solvable group with a partition
\eqref{eq:partition} whose indices are distinct, let $p\geq5$ be the
largest prime divisor of $|G|$, and let $F$ be the subgroup generated
by the $p$-elements of $G$; use the notation of
Definition~\ref{def:local} for this $F$. If some part is not universal
relative to $F$, then
\begin{equation}\label{eq:criterion}
 1\leq\Pi_p\left(S_p+\frac1{p-1}\right).
\end{equation}
\end{proposition}
\begin{proof}
Take $\mathcal P$ to be the set of primes at most $p$, so that
$\Pi=\frac p{p-1}\Pi_p$ and $c_p\Pi=\Pi_p\bigl(S_p+\frac1{p-1}\bigr)$.
Write $f_q=\ord_p(q)$. For $i\in J$ the local index $\ell_i$ is the index of a
proper subgroup $H_i\cap F$ of $F$. By Lemmas~\ref{lem:structure}
and~\ref{lem:cyclotomic}, either $p\mid\ell_i$ or
$q^{f_q}\mid\ell_i$ for a prime $q\ne p$ dividing $|G|$, and such a $q$ is
less than $p$. So every vector $v(d_i)$ with $i\in J$ lies in
\[
 (U+e_p)\ \cup\ \bigcup_{q<p}(U+f_qe_q).
\]
The vectors $v(d_i)$, $i\in J$, are distinct, and their total weight is
$\nu$. For each
$q<p$, the intersection of $U+e_p$ and $U+f_qe_q$ contains
$U+e_p+f_qe_q$, since $U$ is upward closed. Counting $U+e_p$ first and
then the portions of the other sets outside it, and using
\eqref{eq:translate}, gives
\begin{align*}
 \nu&\leq W(U+e_p)+
       \sum_{q<p}W\bigl((U+f_qe_q)\setminus(U+e_p)\bigr)\\
 &\leq\left[\frac1p+
       \sum_{q<p}\left(q^{-f_q}-p^{-1}q^{-f_q}\right)\right]W(U)
 =c_pW(U).
\end{align*}
By Lemma~\ref{lem:shadow}, which applies Sun's bound in $G/F$, the right
side is at most
$c_p\Pi\nu=\Pi_p\bigl(S_p+\frac1{p-1}\bigr)\nu$. As $\nu>0$, the
assertion follows.
\end{proof}

\begin{lemma}\label{lem:numeric}
For every prime $p\geq5$,
\[
 \Pi_p\left(S_p+\frac1{p-1}\right)<1
 \quad\Longleftrightarrow\quad p\ne7.
\]
\end{lemma}
\begin{proof}
At $5$ the expression is $421/432$; at $7$ it is
\[
 \frac{15}{4}\left(\frac18+\frac1{729}+\frac1{15625}+\frac16\right)
 \approx1.099>1.
\]
For the other primes up to $31$ the following upper bounds suffice. Each
bound on $S_p$ is obtained by finding the least $f\geq1$ with
$q^f\equiv1\pmod p$ for each $q<p$ and adding $q^{-f}$. For example,
modulo $13$ the orders of $2,3,5,7,11$ are $12,3,4,12,12$, so
$S_{13}=3^{-3}+5^{-4}+2^{-12}+7^{-12}+11^{-12}<389/10000$; modulo $29$
the two largest terms are $17^{-4}$ and $7^{-7}$ (the orders of $17$ and
$7$ are $4$ and $7$), and the other terms are below $10^{-8}$ in total,
giving $S_{29}<1/70000$.
\begin{center}
\begin{tabular}{ccl}
\toprule
$p$ & upper bound for $\Pi_p$ & upper bound for $S_p$\\
\midrule
$11$ & $35/8$ & $1/180$\\
$13$ & $77/16$ & $389/10000$\\
$17$ & $21/4$ & $1/250$\\
$19$ & $28/5$ & $1/270$\\
$23$ & $59/10$ & $1/2000$\\
$29$ & $31/5$ & $1/70000$\\
$31$ & $317/50$ & $1/32+1/125+1/7000$\\
\bottomrule
\end{tabular}
\end{center}

Let $q<p$ be prime, with $p\geq5$. Since $q^{\ord_p(q)}\equiv1\pmod p$
and $q^{\ord_p(q)}>1$, we have $q^{\ord_p(q)}\geq p+1$. Moreover
$\ord_p(q)\geq3$: order $1$ would give $p\mid q-1$, and order $2$ would
give $p\mid q+1$, hence $q=p-1$, which is even and at least $4$. So
\begin{equation}\label{eq:orderfree}
 S_p\leq\sum_{q<p}\min\{q^{-3},(p+1)^{-1}\}.
\end{equation}
For $37\leq p\leq251$, monotonicity of $\Pi_p$ and
\eqref{eq:orderfree} bound the expression in the lemma by
\begin{equation}\label{eq:range}
 \Pi_{251}\left(\sum_{q<251}\min\{q^{-3},38^{-1}\}+\frac1{36}\right)
 <\frac{9244}{10000},
\end{equation}
by summing the rational terms. No prime lies strictly between $251$ and
$256$.

\smallskip\noindent\emph{The range $p\geq256$.}
We first bound $\Pi_p$. Every prime greater than $7$ is of the form
$6m-1$ or $6m+1$ with $m\geq2$. Put $L=\lceil p/6\rceil$.
All primes $7<q<p$ occur among these pairs for $2\leq m\leq L$.
Including both entries of every pair, even when an entry is composite
or at least $p$, only increases the product, since each factor is
greater than $1$. Hence
\[
 \Pi_p\leq\frac{35}{8}
 \prod_{m=2}^{L}\frac{(6m-1)(6m+1)}{(6m-2)6m}.
\]
For $m\geq2$ the pair factor satisfies
\[
 \left(\frac{(6m-1)(6m+1)}{(6m-2)6m}\right)^3
 \leq\frac{12m+5}{12m-7}.
\]
To verify this inequality, write the right side minus the left side
over the positive denominator $(12m-7)((6m-2)6m)^3$.
Its numerator is
\[
 (12m+5)((6m-2)6m)^3
 -(12m-7)((6m-1)(6m+1))^3.
\]
Writing $M=m-1\geq1$, expansion gives
\[
 29808M^4+109296M^3+149796M^2+90948M+20633>0.
\]
The resulting product telescopes, because $12m-7=12(m-1)+5$:
\[
 \prod_{m=2}^{L}\frac{12m+5}{12m-7}
 =\frac{29}{17}\cdot\frac{41}{29}\cdots
   \frac{12L+5}{12L-7}
 =\frac{12L+5}{17}.
\]
Since $L\leq p/6+1$, we obtain
\begin{equation}\label{eq:wheelprod}
 \Pi_p^3\leq\left(\frac{35}8\right)^3\frac{12L+5}{17}
 \leq\left(\frac{35}8\right)^3\frac{2p+17}{17}.
\end{equation}

Next let $t$ be the integer with $t^3\leq p<(t+1)^3$.
As $p\geq256>6^3$, we have $t\geq6$.
In \eqref{eq:orderfree}, bound each of the at most $t$ terms with
$q\leq t$ by $(p+1)^{-1}$, and the remaining terms by $q^{-3}$.
Extending the latter sum to all integers greater than $t$ gives
\[
 S_p\leq\frac{t}{p+1}+\sum_{n=t+1}^{\infty}n^{-3}.
\]
For $n\geq2$,
\[
 \frac12\bigl((n-1)^{-2}-n^{-2}\bigr)-n^{-3}
 =\frac{3n-2}{2n^3(n-1)^2}>0.
\]
Consequently, for $N>t$,
\[
 \sum_{n=t+1}^{N}n^{-3}
 \leq\frac12\sum_{n=t+1}^{N}\bigl((n-1)^{-2}-n^{-2}\bigr)
 =\frac1{2t^2}-\frac1{2N^2}.
\]
Letting $N\to\infty$ and using $p+1\geq t^3$ yields
\[
 S_p\leq\frac{t}{p+1}+\frac1{2t^2}
 \leq\frac1{t^2}+\frac1{2t^2}.
\]
Also $p-1\geq t^3-1\geq6t^2-1\geq5t^2$, so
\[
 S_p+\frac1{p-1}
 \leq\frac1{t^2}+\frac1{2t^2}+\frac1{5t^2}
 =\frac{17}{10t^2}.
\]
Since $p\geq256$, we have $2p+17\leq(529/256)p$.
Thus \eqref{eq:wheelprod} gives
\[
 \Pi_p^3\leq\left(\frac{35}8\right)^3\frac{529}{4352}p
 <\left(\frac{35}8\right)^3\frac{529}{4352}(t+1)^3.
\]
Finally, $(t+1)/t^2=1/t+1/t^2\leq7/36$ for $t\geq6$, and hence
\begin{align*}
 \left(\Pi_p\left(S_p+\frac1{p-1}\right)\right)^3
 &<\left(\frac{35}8\cdot\frac{17}{10}\cdot
                \frac{t+1}{t^2}\right)^3\frac{529}{4352}\\
 &\leq\left(\frac{833}{576}\right)^3\frac{529}{4352}
 <\left(\frac32\right)^3\frac{529}{4352}<\frac12.
\end{align*}
Taking cube roots proves the required bound for $p\geq256$.
\end{proof}

All the finite inequalities in this proof are also checked by exact
rational arithmetic in the formalization. For the row $p=31$, the
orders $\ord_{31}(q)$ for $q=2,3,5,7,11,13,17,19,23,29$ are
$5,30,3,15,30,30,30,15,10,10$. Remark~\ref{rem:shell} shows that the failure at $p=7$
cannot be removed merely by sharpening the uniform weight bound for the
same divisibility shells.

The Python script \texttt{solvable\_paper\_numerics.py}, provided as an
ancillary file with the arXiv version of this paper, checks the small-prime
values, every bound in the table, the interval estimate \eqref{eq:range},
and the polynomial and telescoping identities in the argument for
$p\geq256$, using exact arithmetic.
It requires only the Python standard library and can be run with
\begin{verbatim}
python solvable_paper_numerics.py
\end{verbatim}

\section{The prime \texorpdfstring{$7$}{7}}\label{sec:seven}

The generic count of Proposition~\ref{prop:criterion} uses only the
divisibility of local indices. At $p=7$ it fails: $8=2^{\ord_7(2)}$ is
one of the relevant cyclotomic divisors, the shifted set $U+v(8)$ has
weight $W(U)/8$, and by Lemma~\ref{lem:numeric} the right side of
\eqref{eq:criterion} is about $1.1$. We therefore use two kinds of
information that the generic count ignores. The first is which minimal
heads (defined below) have no part of local index exactly $7$ or exactly $8$; these give
the missing points. The second is which parts meet the same $F$-coset as
a part of local index $7$ or $8$ (Lemma~\ref{lem:chambers}). Together
these yield incompatible lower and upper bounds involving the weights of
the missing points.

In this section $G$, the partition, $F$, $J$, $U$ and $\nu$ are as in
Proposition~\ref{prop:criterion} with $p=7$, and
$\mathcal P=\{2,3,5,7\}$. In particular,
\[
 r=\ind GF,\qquad d_i=h_i\ell_i,\qquad J=\{i:\ell_i>1\},\qquad
 \nu=\sum_{i\in J}\frac1{d_i},\qquad U=\bigcup_{i\in J}\up v(h_i).
\]
We write
\begin{equation}\label{eq:sevenconstants}
 \Pi=\frac76\Pi_7=\frac{35}{8},\qquad \alpha=e_7=v(7),\qquad
 \beta=3e_2=v(8),\qquad \eps=3^{-6}+5^{-6}.
\end{equation}
Since $\ord_7(2)=3$ and $\ord_7(3)=\ord_7(5)=6$,
Lemma~\ref{lem:cyclotomic} shows that the local index of every
non-universal part is divisible by at least one of $7$, $8$, $3^6$,
$5^6$.

After Lemma~\ref{lem:chambers}, which needs only the notation introduced
so far, we introduce the minimal heads $\mathcal H$, the sets $A$ and $B$,
the missing points $Z_7$ and $Z_8$, and the sets $D$ and $T$;
these are used throughout the rest of the section.

\begin{lemma}\label{lem:chambers}
Let two distinct parts meet the same $F$-coset. Then their local
indices are not coprime. If one of them has local index $7$, the other
has local index divisible by $7$. If one of them has local index $8$,
the other has local index divisible by $7$ or by $8$.
\end{lemma}
\begin{proof}
Both parts are non-universal, since a universal part contains every
$F$-coset it meets. If $g_iH_i$ meets $yF$ in a point $y'$, then
$g_iH_i\cap yF=y'(H_i\cap F)$. So two parts $i,j$ meeting $yF$ give
disjoint pieces $y'(H_i\cap F)$ and $y''(H_j\cap F)$. Translating by
$y'^{-1}$, we obtain disjoint cosets $Y_0$ and $uY$ in $F$, where
$Y_0=H_i\cap F$, $Y=H_j\cap F$, and $u=(y')^{-1}y''\in F$ because
$y'$ and $y''$ lie in the same left coset $yF$. The indices are the local
indices $\ell_i,\ell_j$. By Lemma~\ref{lem:coprime}, $\ell_i$ and
$\ell_j$ are not coprime. This proves the first assertion, and the
second follows because $7$ is prime.

For the third, suppose $\ind F{Y_0}=8$. If $7\mid\ind FY$ there is nothing
to prove. Otherwise $Y$ contains a Sylow $7$-subgroup $P$ of $F$. The
action of $P$ on the eight cosets of $Y_0$ is nontrivial. Indeed, if $P$
lay in the core of $Y_0$, that normal subgroup of $F$ would have index
prime to $7$, so it would be $F$ by Lemma~\ref{lem:structure}, contrary to
$\ind F{Y_0}=8$.
A nontrivial $7$-group action on eight points has one orbit of size $7$
and one fixed point $x_0$. Since $P\leq Y$, the $Y$-orbit of $x_0$ and
its complement are $P$-invariant, so the complement has size $0$ or $7$.
In the first case $Y$ is transitive on $F/Y_0$, so $F=YY_0$. Then
$uY\cap Y_0\ne\varnothing$, since $u^{-1}\in F=YY_0$, contrary to the two
disjoint pieces. In the second case $Y$ fixes $x_0$ and is contained in
its stabilizer, a conjugate of $Y_0$. Thus $8\mid\ind FY$.
\end{proof}

We call $v(h_i)$ the \emph{head} of part $i$. Let $\mathcal H$ be the
set of minimal elements of $\{v(h_i):i\in J\}$, so that $\mathcal H$ is
a finite antichain and $U=\up\mathcal H$. We call the elements of
$\mathcal H$ the \emph{minimal heads}. Let $A\subseteq\mathcal H$
consist of those $\eta$ for which some part has head $\eta$ and local
index exactly $7$, and let $B\subseteq\mathcal H$ consist of those
$\eta$ for which some part has head $\eta$ and local index exactly $8$.
Set
\[
 Z_7=\{\eta+\alpha:\eta\in\mathcal H\setminus A\},\qquad
 Z_8=\{\eta+\beta:\eta\in\mathcal H\setminus B\},\qquad
 z_7=W(Z_7),\quad z_8=W(Z_8).
\]
We call the elements of $Z_7\cup Z_8$ the \emph{missing points}.
A missing point may still be the
vector of some other part if it lies in both $U+\alpha$ and $U+\beta$.
To handle this we use the overlap set
\[
 D=((U+\alpha)\cap(U+\beta))\setminus(U+\alpha+\beta),
\]
and we put
\[
 T=(U+\alpha)\cup(U+\beta).
\]
Inclusion--exclusion and \eqref{eq:translate} give
\begin{equation}\label{eq:defect}
 W(T)+W(D)=\left(\frac17+\frac18-\frac1{56}\right)W(U)=\frac14W(U).
\end{equation}

\begin{lemma}[Missing points in the overlap]\label{lem:charging}
With the notation above:
\begin{enumerate}
\item $\displaystyle\frac74 W(Z_7\cap D)+W((Z_8\setminus Z_7)\cap D)\leq W(D)$;
\item $Z_7\cap Z_8\subseteq D$;
\item if $i\in J$ and $7\mid\ell_i$ or $8\mid\ell_i$, then $v(d_i)$ lies
in neither $Z_7\setminus D$ nor $Z_8\setminus D$.
\end{enumerate}
\end{lemma}
\begin{proof}
We use repeatedly that for $\eta\in\mathcal H$ and $x\in U$,
$x\leq\eta$ implies $x=\eta$: some $\eta'\in\mathcal H$ has
$\eta'\leq x\leq\eta$, and $\eta'=\eta$ because $\mathcal H$ is an
antichain. In particular, for $\eta\in\mathcal H$ neither $\eta+\alpha$
nor $\eta+\beta$ lies in $U+\alpha+\beta$, since
$\eta+\alpha=x+\alpha+\beta$ or $\eta+\beta=x+\alpha+\beta$ with
$x\in U$ would give $x<\eta$.

(1) For each $\eta+\alpha\in Z_7\cap D$ and $m\in\{0,1,2\}$, consider
$c_m(\eta)=\eta+\alpha+me_2$. Each $c_m(\eta)$ lies in $U+\alpha$,
because $\eta+me_2\in U$, and in $U+\beta$, because
$\eta+\alpha\in D\subseteq U+\beta$ and $U+\beta$ is upward closed. None
lies in $U+\alpha+\beta=U+\alpha+3e_2$: an equality
$c_m(\eta)=x+\alpha+3e_2$ with $x\in U$ gives $x=\eta-(3-m)e_2<\eta$.
Hence $c_0(\eta),c_1(\eta),c_2(\eta)\in D$, and their weights add to
$(1+1/2+1/4)w(\eta+\alpha)=\frac74w(\eta+\alpha)$.

These triples are disjoint. If $c_m(\eta)=c_{m'}(\eta')$ then
$\eta-\eta'=(m'-m)e_2$, so $\eta,\eta'\in\mathcal H$ are comparable; the
antichain property gives $\eta=\eta'$ and then $m=m'$. They also avoid
$(Z_8\setminus Z_7)\cap D$. The point $c_0(\eta)=\eta+\alpha$ lies in
$Z_7$. For $m\in\{1,2\}$ suppose $c_m(\eta)=\eta'+\beta$ with
$\eta'\in\mathcal H$, so $\eta'=\eta+\alpha-(3-m)e_2$. Since
$\eta+\alpha\in U+\beta$, some $\eta''\in\mathcal H$ satisfies
$\eta''\leq\eta+\alpha-3e_2<\eta'$, contradicting the antichain
property. The disjoint triples and $(Z_8\setminus Z_7)\cap D$ all lie in
$D$, which proves (1).

(2) A point of $Z_7\cap Z_8$ lies in $(U+\alpha)\cap(U+\beta)$, and by
the first paragraph it does not lie in $U+\alpha+\beta$.

(3) Put $x=v(h_i)\in U$ and $l=v(\ell_i)$, so that $v(d_i)=x+l$.
Suppose first that $v(d_i)\in Z_7\setminus D$, so that
$x+l=\eta+\alpha$ for some $\eta\in\mathcal H\setminus A$.
If $7\mid\ell_i$, then $l\geq\alpha$, hence $x\leq\eta$.
Thus $x=\eta$ and $l=\alpha$, so part $i$ has head $\eta$ and local
index $7$, contrary to $\eta\notin A$.
If $8\mid\ell_i$, then $l\geq\beta$, so $x+l\in U+\beta$.
It also lies in $U+\alpha$, and by the first paragraph it does not lie
in $U+\alpha+\beta$. Hence $x+l\in D$, again a contradiction.

Now suppose that $v(d_i)\in Z_8\setminus D$, so that
$x+l=\eta+\beta$ for some $\eta\in\mathcal H\setminus B$.
If $8\mid\ell_i$, then $l\geq\beta$, hence $x\leq\eta$.
Thus $x=\eta$ and $l=\beta$, so part $i$ has head $\eta$ and local
index $8$, contrary to $\eta\notin B$.
If $7\mid\ell_i$, then $l\geq\alpha$, so $x+l\in U+\alpha$.
It also lies in $U+\beta$, and by the first paragraph it does not lie
in $U+\alpha+\beta$. Hence $x+l\in D$, a contradiction.
\end{proof}

\begin{proposition}\label{prop:seven}
Let $G$ be a finite solvable group with a partition
\eqref{eq:partition} whose indices are distinct, suppose that $7$ is the
largest prime divisor of $|G|$, and let $F$ be the subgroup generated
by the $7$-elements of $G$. Then every part
is universal relative to $F$.
\end{proposition}
\begin{proof}
Suppose, for a contradiction, that $J\ne\varnothing$. Then
$U\ne\varnothing$. Put
\[
 \lambda=\frac{W(U)}\Pi>0,\qquad
 z_{7,\mathrm{out}}=W(Z_7\setminus D),\qquad
 z_{8,\mathrm{out}}=W(Z_8\setminus D).
\]
By \eqref{eq:shadow}, $\lambda\leq\nu$.

\smallskip\noindent\emph{Step 1: a lower bound for the missing points.}
Let $\Omega_7\subseteq G/F$ be the union of the images $X_i$ of the parts
with local index exactly $7$ and head in $A$, and define $\Omega_8$
likewise using local index $8$ and $B$. Put $\rho_7=|\Omega_7|/r$ and
$\rho_8=|\Omega_8|/r$. Membership of a head in $A$ or $B$ asserts
the existence of a part of local index $7$ or $8$ with that head;
we select only parts with the specified local index.
By Lemma~\ref{lem:shadow}, $W(\up A)\leq\Pi\rho_7$
and $W(\up B)\leq\Pi\rho_8$. Parts of local index $7$ or $8$ whose heads
are not minimal are not needed here. Since
$U=\up A\cup\bigcup_{\eta\in\mathcal H\setminus A}\up\eta$ and
$W(\up\eta)=\Pi w(\eta)=7\Pi w(\eta+\alpha)$ by \eqref{eq:translate},
\[
 W(U)\leq\Pi\rho_7+7\Pi\sum_{\eta\in\mathcal H\setminus A}w(\eta+\alpha)
 =\Pi(\rho_7+7z_7).
\]
With $w(\eta)=8w(\eta+\beta)$ the same argument for $B$ gives
\begin{equation}\label{eq:apexlower}
 \lambda\leq\rho_7+7z_7,\qquad \lambda\leq\rho_8+8z_8.
\end{equation}

The sets $\Omega_7$ and $\Omega_8$ are disjoint: an $F$-coset in both
would be met by parts of local indices $7$ and $8$, contrary to
Lemma~\ref{lem:chambers}. Let $I_\Omega$ be the set of parts meeting the
inverse image of $\Omega_7\cup\Omega_8$ in $G$, and let
$V_\Omega=\{v(d_i):i\in I_\Omega\}$. These parts are non-universal,
because each such $F$-coset is met by a non-universal part. Their union
contains the inverse image of $\Omega_7\cup\Omega_8$.
For the quotient map $\pi:G\to G/F$, this inverse image has proportion
\[
 \frac{|\pi^{-1}(\Omega_7\cup\Omega_8)|}{|G|}
 =\frac{|F|(|\Omega_7|+|\Omega_8|)}{|F|r}
 =\rho_7+\rho_8.
\]
As the $d_i$ are distinct,
\[
 \rho_7+\rho_8\leq\sum_{i\in I_\Omega}\frac1{d_i}=W(V_\Omega).
\]
By Lemma~\ref{lem:chambers}, the local index of each part in $I_\Omega$
is divisible by $7$ or $8$. Hence $V_\Omega\subseteq T$, and
Lemma~\ref{lem:charging}(3) shows that $V_\Omega$ does not meet
$Z_7\setminus D\subseteq T$. Therefore, using \eqref{eq:defect},
\begin{equation}\label{eq:K3}
 \rho_7+\rho_8+z_{7,\mathrm{out}}\leq W(T)=\frac14W(U)-W(D)
 =\frac{35}{32}\lambda-W(D).
\end{equation}
Adding the two inequalities \eqref{eq:apexlower} and using
\eqref{eq:K3},
\[
 2\lambda\leq\rho_7+\rho_8+7z_7+8z_8
 \leq\frac{35}{32}\lambda-W(D)-z_{7,\mathrm{out}}+7z_7+8z_8.
\]
Since $z_7=W(Z_7\cap D)+z_{7,\mathrm{out}}\leq W(D)+z_{7,\mathrm{out}}$,
this gives $\frac{29}{32}\lambda\leq6z_7+8z_8$, that is,
\begin{equation}\label{eq:forced}
 \frac34z_7+z_8\geq\frac{29}{256}\lambda.
\end{equation}

\smallskip\noindent\emph{Step 2: an upper bound for the missing points.}
Divide the non-universal parts by local divisibility:
\[
 J_{78}=\{i\in J:7\mid\ell_i\text{ or }8\mid\ell_i\},\qquad
 J_{35}=J\setminus J_{78},
\]
and let $V_{78}$ and $V_{35}$ be the corresponding sets of vectors
$v(d_i)$. They are disjoint because the $d_i$ are distinct, and
$\nu=W(V_{78})+W(V_{35})$. We have $V_{78}\subseteq T$, and by
Lemma~\ref{lem:charging}(3) the set $V_{78}$ meets neither
$Z_7\setminus D$ nor $Z_8\setminus D$. These two sets lie in $T$ and are
disjoint by Lemma~\ref{lem:charging}(2). Hence
$W(V_{78})+z_{7,\mathrm{out}}+z_{8,\mathrm{out}}\leq W(T)$. For
$i\in J_{35}$, the local index is divisible by $3^6$ or $5^6$, so
$V_{35}\subseteq(U+6e_3)\cup(U+6e_5)$ and $W(V_{35})\leq\eps W(U)$.
Adding and using \eqref{eq:defect},
\[
 \nu+z_{7,\mathrm{out}}+z_{8,\mathrm{out}}\leq W(T)+\eps W(U)
 =\frac14W(U)-W(D)+\eps W(U).
\]
Since $\nu\geq\lambda$ by Lemma~\ref{lem:shadow} and
$W(U)=\frac{35}8\lambda$, this gives
\begin{equation}\label{eq:K4}
 W(D)+z_{7,\mathrm{out}}+z_{8,\mathrm{out}}\leq
 \left(\frac3{32}+\frac{35}{8}\eps\right)\lambda.
\end{equation}

The coefficient $\frac34=\frac68$ in \eqref{eq:forced} comes from
$7z_7-z_7$ and $8z_8$. A point of $Z_7\cap D$ may also lie in $Z_8$, so
below it contributes up to $\frac34+1=\frac74$ times its weight, which
is exactly the factor in Lemma~\ref{lem:charging}(1).

\smallskip\noindent\emph{Step 3: conclusion.}
Put $\gamma=W(Z_7\cap D)$ and $\delta=W((Z_8\setminus Z_7)\cap D)$.
Then $z_7=z_{7,\mathrm{out}}+\gamma$. By Lemma~\ref{lem:charging}(2),
$Z_7\cap Z_8\subseteq Z_7\cap D$, so splitting $Z_8$ into its parts
outside $D$, in $D\setminus Z_7$, and in $Z_7$ gives
$z_8\leq z_{8,\mathrm{out}}+\delta+\gamma$. By
Lemma~\ref{lem:charging}(1),
\[
 \frac34z_7+z_8
 \leq\frac34z_{7,\mathrm{out}}+z_{8,\mathrm{out}}+\frac74\gamma+\delta
 \leq\frac34z_{7,\mathrm{out}}+z_{8,\mathrm{out}}+W(D)
 \leq z_{7,\mathrm{out}}+z_{8,\mathrm{out}}+W(D).
\]
Combining this with \eqref{eq:forced} and \eqref{eq:K4} gives
\[
 \frac{29}{256}\leq\frac3{32}
       +\frac{35}{8}\left(\frac1{729}+\frac1{15625}\right),
\]
which is false: the left side exceeds the right side by
$7727329/583200000$.
\end{proof}

\begin{remark}\label{rem:shell}
The shortfall is not an artifact of the union bound in the proof of
Proposition~\ref{prop:criterion}. Let $\mathcal P=\{2,3,5,7\}$, list the
shifts $\alpha,\beta,6e_3,6e_5$ as $a_1,\ldots,a_4$, and let $U\subseteq\N^{\mathcal P}$
be nonempty and upward closed. Put $V_m=U\setminus\bigcup_{j<m}(U+a_j)$.
Since $U$ is upward closed, $V_m\cap(U+a_m)\subseteq V_m+a_m$, so
$W(V_{m+1})\geq(1-w(a_m))W(V_m)$. Iterating, and using $U+a_j\subseteq U$,
\[
 W\Bigl(\bigcup_m(U+a_m)\Bigr)\leq c\,W(U),\qquad
 c=1-\prod_m\bigl(1-w(a_m)\bigr).
\]
For $U=\N^{\mathcal P}$ the four shifts lie on distinct coordinates, so
$U\setminus\bigcup_m(U+a_m)$ is a product and equality holds. Thus $c$ is
the best uniform constant, and still $\frac{35}8c\approx1.098>1$.
\end{remark}

\section{Proof of the main theorem}\label{sec:finish}

\begin{proof}[Proof of Theorem~\ref{thm:main}]
By Lemma~\ref{lem:finite} it suffices to treat finite groups. Suppose that
$G$ is a finite solvable group of least order that is not HS, and fix a
partition with distinct indices. By Lemma~\ref{lem:mincex}, the largest
prime divisor $p$ of $|G|$ is at least $5$. Let $F$ be the subgroup
generated by the $p$-elements of $G$. By Lemma~\ref{lem:structure}, it
is a nontrivial normal subgroup. If every part is universal relative to
$F$, the partition descends to a partition of the smaller solvable group
$G/F$ with the same indices, contradicting minimality. Otherwise
Proposition~\ref{prop:criterion} and Lemma~\ref{lem:numeric} exclude
$p\ne7$, and Proposition~\ref{prop:seven} excludes $p=7$.
\end{proof}

\begin{remark}\label{rem:solvable}
Solvability is used in two places: in Theorem~\ref{thm:sun}, applied to
the quotient $G/F$, and in Lemma~\ref{lem:cyclotomic}, where a minimal
normal subgroup of the primitive quotient $F/R$ is elementary abelian, so
that the degree $\ind FM$ is a prime power. Lemma~\ref{lem:structure}
holds for every finite group. The counting in Sections~\ref{sec:general}
and~\ref{sec:seven} uses only the conclusions of these results. Sun's
theorem also holds without solvability for cosets of subnormal subgroups
\cite[Theorem~3.1]{Sun}, but the images $H_iF/F$ need not be subnormal,
and for nonsolvable groups primitive degrees need not be prime powers, so
Lemma~\ref{lem:cyclotomic} would also need replacement.
\end{remark}

\begin{remark}\label{rem:localHS}
The only direct uses of minimality of $G$ are to exclude index $2$ and
to descend to $G/F$. In the setting of the proof minimality gives more when $F<G$.
The parts meeting a non-universal $F$-coset $yF$ number at least two, and
after translation by an element of $yF$ their pieces
$y'(H_i\cap F)$ form a partition of $F$ into left cosets. As $F$ is a
smaller solvable group, it is HS, so there exist two distinct parts
meeting $yF$ with equal local index. We do not use this, because the
proof cannot
assume $F<G$, including when $p=7$. For example, in $C_2^3\rtimes C_7$ the
subgroup generated by the $7$-elements is the whole group when $C_7$ acts
nontrivially and irreducibly on $C_2^3$. This semidirect product is unique
up to isomorphism, since $\operatorname{GL}(3,2)$ has a unique conjugacy
class of cyclic subgroups of order $7$. This group is HS, so it is not a
counterexample, but
largest prime $7$ and solvability do not by themselves force $F<G$.
\end{remark}

\section{Formal proof}\label{sec:formal}

The accompanying Lean development \cite{Repository} proves
Theorem~\ref{thm:main} as
\begin{center}
\texttt{Erdos274.herzog\_schonheim.variants.solvable}.
\end{center}
The proof revision is
\href{https://github.com/Jostamon/erdos274-hs-solvable/commit/d701decf37c225a3e5326fade4bd92121b1f93d4}{\texttt{d701decf37c2}}.
It is pinned to Lean \texttt{v4.35.0-rc2} and the mathlib tag
\texttt{v4.35.0-rc2}, a release candidate, by the files
\texttt{lean-toolchain} and \texttt{lakefile.toml}. The repository consists
of exactly the import closure of the main theorem. From its root, build it
with
\begin{verbatim}
lake exe cache get
lake build
\end{verbatim}
In a Lean file importing the module \texttt{Erdos274.KnownCases.Solvable},
the command
\begin{verbatim}
#print axioms Erdos274.herzog_schonheim.variants.solvable
\end{verbatim}
reports precisely
\begin{verbatim}
[propext, Classical.choice, Quot.sound]
\end{verbatim}
Thus the elaborated theorem has no further axiomatic dependencies,
including \texttt{sorryAx}. This check concerns the formal theorem;
the interpretation of its hypotheses and its correspondence with the
paper remain mathematical questions. In mathematical notation, the
formal theorem reads as follows. Let $G$ be a solvable group with more
than one element, and let $\iota$ be a finite type with at least two
elements. Let $(H_i)_{i\in\iota}$ be subgroups and $(g_i)_{i\in\iota}$
elements of $G$ such that the left cosets $g_iH_i$ are nonempty, pairwise
disjoint, and have union $G$. Then $\ind G{H_i}=\ind G{H_j}$ for some
$i\ne j$. No finiteness hypothesis is placed on $G$, and the cosets are
left cosets, as in \eqref{eq:partition}. The index is mathlib's
\texttt{Subgroup.index}, which is $0$ for a subgroup of infinite index;
in a finite partition every index is finite (Lemma~\ref{lem:finite}).
The hypothesis $|G|>1$ follows from the others and is kept to match the
upstream statement. The exact covering structure is copied field for field
from the Erd\H{o}s problem~274 file of the Formal Conjectures repository
\cite{FormalConjectures,FormalConjecturesPaper},
and the file \texttt{Challenge.lean} of the repository states the theorem
using only mathlib and that structure, so the formal statement can be
audited without reading the proof.

The formalization follows Sections~\ref{sec:local}--\ref{sec:seven}.
Its treatment of the inputs from
Section~\ref{sec:prelim} differs in form only. Lean proves Sun's bound itself, in a more general form for
orbits, and states it as
$\sum_{d\in\mathcal D}\varphi(d)\leq\bigl|\bigcup_jx_jL_j\bigr|$, where
$\mathcal D$ is the set of divisors of the orders $|L_j|$. Counting
elements of a cyclic group of order $|L|$ by their orders shows that the
left side equals the left side of Theorem~\ref{thm:sun}. Lemma~\ref{lem:dictionary}
is proved there by decomposing $U$ into fibres over the divisors of a
common multiple, rather than by inclusion--exclusion. Likewise, to see that
every part of a partition of a possibly infinite group has finite index,
Lean uses mathlib's version of B.~H.~Neumann's lemma
\cite[Lemma~4.1]{Neumann}, in place of \cite[Theorem~1]{KZ}.

\section*{Acknowledgments}
\phantomsection\label{sec:ack}

The author thanks Leo Margolis for a careful reading of an earlier version,
many helpful comments, and his encouragement.

Under the author's direction, OpenAI GPT-5.6 Sol contributed the
comparison of reciprocal sums with unions of quotient cosets
(Section~\ref{sec:local}), the cyclotomic divisibility argument, and
the argument for the prime $7$, including the overlap set $D$ and
Lemma~\ref{lem:charging}. Anthropic Claude Opus 5 developed
and formalized the early solvable track and the enlargement to all
non-universal parts. Claude Opus 5.5 audited and formalized the final
argument and completed the theorem for arbitrary solvable groups; its
simplifications include localizing at the subgroup generated by the
$p$-elements, the fixed-point proofs of the local
divisibility restrictions, including the maximal-subgroup proof of
Lemma~\ref{lem:cyclotomic}, the bound $\ord_p(q)\geq3$, and the
elementary estimates in Section~\ref{sec:general}.
OpenAI GPT-5.6 Sol, reviewing the paper, proposed the present form of the
large-prime estimate in the proof of Lemma~\ref{lem:numeric}. Claude Fable 5
contributed earlier known-case formalizations and conditional
reductions. Reviews by DeepSeek (version not recorded) and Z.ai GLM 5.2
clarified the earlier cyclotomic argument. OpenAI GPT-6 in Codex assisted
with the exposition and attribution audit. The commit history of the
repository \cite{Repository} records the AI contributions to the
formalization through its \texttt{Co-Authored-By} trailers.

The author thanks the Lean and mathlib communities for the libraries
supporting the formal proof. The author is responsible for the
mathematical claims, their attribution, and the correspondence between
the paper and the formal statements.


\begin{thebibliography}{99}

\bibitem{AS25a}
F. Akman and P. Sissokho,
\emph{Transversal coset partitions of groups},
Beitr. Algebra Geom. \textbf{66} (2025), 417--441,
\href{https://doi.org/10.1007/s13366-024-00748-9}{doi:10.1007/s13366-024-00748-9}.

\bibitem{AS25b}
F. Akman and P. Sissokho,
\emph{Steiner coset partitions of groups},
Canad. Math. Bull. \textbf{68} (2025), 1359--1373,
\href{https://doi.org/10.4153/S0008439525100787}{doi:10.4153/S0008439525100787}.

\bibitem{AS26}
F. Akman and P. Sissokho,
\emph{Steiner coset partitions for five mutually commuting subgroups},
Acta Math. Hungar. (2026), published online,
\href{https://doi.org/10.1007/s10474-026-01635-6}{doi:10.1007/s10474-026-01635-6}.

\bibitem{BFF86}
M.~A. Berger, A. Felzenbaum, and A.~S. Fraenkel,
\emph{The Herzog--Sch\"onheim conjecture for finite nilpotent groups},
Canad. Math. Bull. \textbf{29} (1986), no.~3, 329--333.

\bibitem{BFF87}
M.~A. Berger, A. Felzenbaum, and A.~S. Fraenkel,
\emph{Remark on the multiplicity of a partition of a group into cosets},
Fund. Math. \textbf{128} (1987), 139--144.

\bibitem{Burkhart}
M.~C. Burkhart,
\emph{The Herzog--Sch\"onheim conjecture for solvable groups},
\href{https://arxiv.org/abs/1901.10131}{arXiv:1901.10131} (2019),
withdrawn in version~2, March~1, 2019.

\bibitem{FormalConjectures}
The Formal Conjectures Authors, \emph{The Formal Conjectures Repository},
2025, \url{https://github.com/google-deepmind/formal-conjectures}.

\bibitem{FormalConjecturesPaper}
M. Firsching, P. Lezeau, S. Mercuri, M.~Z. Horv\'ath, Y. Dillies,
C. S\"onne, E. Wieser, F. Zhang, T. Hubert, B. Ag\"uera y Arcas, and
P. Kohli, \emph{Formal Conjectures: An open and evolving benchmark for
verified discovery in mathematics}, 2026,
\href{https://arxiv.org/abs/2605.13171v1}{arXiv:2605.13171v1}.

\bibitem{GM}
M. Garonzi and L. Margolis,
\emph{The Herzog--Sch\"onheim conjecture for simple and symmetric groups},
Algebra \& Number Theory, to appear;
\href{https://arxiv.org/abs/2509.25118v2}{arXiv:2509.25118v2}.

\bibitem{GS}
Y. Ginosar and O. Schnabel,
\emph{Prime factorization conditions providing multiplicities in coset
partitions of groups}, J. Comb. Number Theory \textbf{3} (2011), no.~2, 75--86.

\bibitem{HS}
M. Herzog and J. Sch\"onheim,
\emph{Research problem no.~9}, Canad. Math. Bull. \textbf{17} (1974), 150.

\bibitem{KZ}
I. Korec and \v S. Zn\'am,
\emph{On disjoint covering of groups by their cosets},
Math. Slovaca \textbf{27} (1977), no.~1, 3--7.

\bibitem{MS}
L. Margolis and O. Schnabel,
\emph{The Herzog--Sch\"onheim conjecture for small groups and harmonic
subgroups}, Beitr. Algebra Geom. \textbf{60} (2019), no.~3, 399--418;
\href{https://arxiv.org/abs/1803.03569}{arXiv:1803.03569}.

\bibitem{Repository}
M. Menon, \emph{The Herzog--Sch\"onheim conjecture for solvable groups in
Lean~4}, source repository,
\url{https://github.com/Jostamon/erdos274-hs-solvable},
commit \texttt{d701decf37c225a3e5326fade4bd92121b1f93d4} (October~2, 2026).

\bibitem{Neumann}
B.~H. Neumann, \emph{Groups covered by permutable subsets},
J. London Math. Soc. \textbf{29} (1954), 236--248,
\href{https://doi.org/10.1112/jlms/s1-29.2.236}{doi:10.1112/jlms/s1-29.2.236}.

\bibitem{Robinson}
D.~J.~S. Robinson, \emph{A Course in the Theory of Groups}, second ed.,
Graduate Texts in Mathematics, vol.~80, Springer, New York, 1996.

\bibitem{Sun06}
Z.-W. Sun, \emph{Finite covers of groups by cosets or subgroups},
Internat. J. Math. \textbf{17} (2006), no.~9, 1047--1064,
\href{https://arxiv.org/abs/math/0501451}{arXiv:math/0501451}.

\bibitem{Sun}
Z.-W. Sun, \emph{On the Herzog--Sch\"onheim conjecture for uniform covers
of groups}, J. Algebra \textbf{273} (2004), 153--175,
\href{https://arxiv.org/abs/math/0306099}{arXiv:math/0306099}.

\end{thebibliography}
\end{document}